\documentclass[12pt]{amsart}
\usepackage{amssymb,amsmath}
\usepackage[T1]{fontenc}
\usepackage[utf8]{inputenc}
\usepackage{enumitem}
\usepackage{hyperref}
\usepackage{xurl}
\usepackage[labelfont=rm]{caption}
\usepackage[labelformat=simple]{subcaption}
\usepackage[defaultlines=3,all]{nowidow}
\usepackage{tikz}
\usepackage{mathrsfs}

\usepackage{pgf}
\usepackage{xcolor}
\usepackage{nicematrix}

\usetikzlibrary{braids, patterns,
shapes, matrix, arrows, knots, arrows.meta, 3d, pgfplots.groupplots, angles, ext.transformations.mirror,
calc, decorations.markings, intersections, cd, external}

\DeclareMathOperator{\chr}{Char}
\DeclareMathOperator{\ord}{ord}
\newcommand{\nd}{{\rm ND}}

\newcommand{\bn}{\mathbb{N}}
\newcommand{\bk}{\mathbb{K}}

\newcommand{\cT}{\mathcal{T}}
\newcommand{\cV}{\mathcal{V}}
\newcommand{\cA}{\mathcal{A}}
\newcommand{\cE}{\mathcal{E}}

\newcommand{\area}{\mathcal{A}}
\newcommand{\arista}{\mathscr{A}}
\newtheorem{theorem}{Theorem}[section]
\newtheorem{proposition}[theorem]{Proposition}

\newtheorem{corollary}[theorem]{Corollary}
\newtheorem{lemma}[theorem]{Lemma}
\newtheorem*{cjt}{Conjecture}

\theoremstyle{definition}
\newtheorem{definition}[theorem]{Definition}
\newtheorem{ntc}[theorem]{Notation}
\newtheorem{example}[theorem]{Example}

\theoremstyle{remark}
\newtheorem{remark}[theorem]{Remark}

\numberwithin{equation}{section}

\makeatletter
\@namedef{subjclassname@2020}{\textup{2020} Mathematics Subject Classification}
\makeatother

\begin{document}

\baselineskip=17pt

\title[Milnor number]{Milnor number of plane curve singularities in arbitrary characteristic}    \author[E. Artal]{Enrique Artal Bartolo}
\address{Departamento de Matem\'{a}ticas, IUMA \\
Universidad de Zaragoza \\
C.~Pedro Cerbuna 12, 50009, Zaragoza, Spain}
\urladdr{\url{http://riemann.unizar.es/~artal}}
\email{\href{mailto:artal@unizar.es}{artal@unizar.es}}

\author[Pi. Cassou-Nogu\`es]{Pierrette Cassou-Noguès}

\address{Institut de Math\'ematiques de Bordeaux, Universit\'e Bordeaux,
350, Cours de la Lib\'eration, 33405, Talence Cedex 05, FRANCE}
\email{\href{mailto:cassou@math.u-bordeaux.fr}{cassou@math.u-bordeaux.fr}}

\date{}

\begin{abstract}
Reduced power series in two variables with coefficients in a field of characteristic zero satisfy
a well-known formula that relates a codimension related to the normalization of a ring and 
the Jacobian ideal. In the general case Deligne proved that this formula is only an inequality;
Garc{\'i}a Barroso and P{\l}oski stated a conjecture for irreducible power series. In this work we 
generalize Kouchnirenko's formula for any reduced power series and also generalize 
Garc{\'i}a Barroso and P{\l}oski's conjecture. We prove the conjecture in some cases using in particular 
Greuel-Nguyen's results.

\end{abstract}

\subjclass[2020]{Primary 14H20, 14G17,14B05, 32A05; Secondary 14M25 14Q05}

\keywords {trees, Newton polygon, field of arbitrary characteristic, $\delta$-invariant, Milnor number}

\maketitle

Let $\bk$ be an algebraically closed field of arbitrary characteristic, and let
$f\in\bk[[x,y]]$ be a reduced power series. Let $\overline{\mathcal{O}}$ be
the normalization of the ring $\mathcal{O}:=\bk[[x,y]]/(f)$,
and let $\delta(f):=\dim_\bk\overline{\mathcal{O}}/\mathcal{O}$.
We set
\[
\overline{\mu}(f):=2\delta(f)-r(f)+1
\]
where $r(f)$ is the number of distinct irreducible factors of $f$.

The main result of this article is the computation of $\overline{\mu}(f)$ in terms of areas of Newton polygons, in the spirit of Kouchnirenko, without any hypothesis on degeneration.
Let
\[
\mu(f):=\dim_\bk \bk[[x,y]]/\left(\frac{\partial f}{\partial x}, \frac{\partial f}{\partial y}\right).
\]
In characteristic zero, we have $\mu(f)=\overline{\mu}(f)$.
Deligne proved that
\[
\mu(f)\geq\overline{\mu}(f).
\]
Some authors as Garc{\'i}a-Barroso and P{\l}oski~\cite{evpl:18}, Greuel and Nguyen \cite{gn:12}, and Hefez, Rodrigues, Solom\~{a}o ~\cite{hrs:18}, were interested in the question of the equality in characteristic $p\neq 0$.  We give a conjecture on this question, and using works of these authors, we prove that after adding some hypothesis, it is true. We show some more examples of its validity.

The paper is organized as follows.
In~\S\ref{sec:hn}, we study the Hamburger-Noether algorithm, defined in a form which is very close to the Newton algorithm, but that can be used in any characteristic.

In~\S\ref{sec:nt}, using the Hamburger-Noether algorithm, we construct trees for any reduced power series in $f(x,y)\in\bk[[x,y]]$ in any characteristic. Note that, given~$f$ the tree depends on $\chr\bk$. These trees are constructed using the Newton polygons of $f$ and 
the ones of its transforms at each stage of the Hamburger-Noether algorithm. We define an important invariant: the multiplicity of the tree. It is defined using the decorations of the tree, which are computed from the equations of the faces of the Newton polygons.

In~\S\ref{sec:area}, we show that in fact, the multiplicity of the tree can be also expressed in terms of the areas below the Newton polygons that appear in the Hamburger-Noether algorithm. This is a generalization of Kouchnirenko's result, see Remark~\ref{rem:kou}.

The following three sections aim to prove the main result of the paper, that is that  the multiplicity of the tree is also equal to $r-2\delta(f)$.  For this we need to compute the multiplicity of intersection of two power series in terms of the trees, which is done 
in~\S\ref{sec:inter}. Then in~\S\ref{sec:irr} we study irreducible series and show how to compute their Zariski characteristic series from the tree. The main result is proven in~\S\ref{sec:delta}, first for irreducible series, and then in general using results from~\cite{cnd}.

The last section~\S\ref{sec:milnor} is devoted to the study of the following conjecture.
Let $f\in\bk[[x,y]]$, and let $\mathcal{T}(f)$ be its minimal tree, see~\S\ref{sec:nt}.
Let $\cV$ be the set of vertices of $\mathcal{T}(f)$.

\begin{cjt}
We have \[
\mu(f)=\overline{\mu}(f)
\]
if and only if $p$ does not divide any of the $N_v$ for all $v\in \mathcal{V}.$
\end{cjt}

The previous sections compute the term on the right hand side.
Using Greuel and Nguyen result~\cite{gn:12}, we show that the conjecture is true when $f$ is non-degenerate. Here one has to be careful about the definition of degenerate, since our definition does not coincide with Greuel-Nguyen's one.

Using a result from Garc{\'i}a-Barroso and P{\l}oski~\cite{evpl:18}, we prove that if $\mathcal{T}$ is a tree and $p>M(\mathcal{T})+\ord (\mathcal{T})$,  for all $f\in\bk[[x,y]]$ with characteristic of $\bk$ equal to $p$ and $\mathcal{T}(f)=\mathcal{T}$, then $\mu(f)=\overline{\mu}(f).$ In this case, $p$ divides~$N_v$ for no vertex~$v$ of $\cT$.

We show also that the conjecture agrees with the conjecture of  Garc{\'i}a-Barroso and P{\l}oski~\cite{evpl:18}, in the case where $f$ is irreducible, and with the result of Hefez, Rodrigues, Solom\~{a}o.
We give some examples where the conjecture is true.

In a subsequent article in preparation, we will show some other parts of the conjecture.

This article owes a great deal to Arkadiusz P{\l}oski for many reasons. The second author met him about thirty years ago in Bordeaux, when he came for a month as invited Professor at the University. It was the beginning of a strong friendship and collaboration. Over the years, she learnt, discussing with him, not only mathematics, but also, history, political sciences, literature,\dots\ She wants also to mention the important part played by his wife, Anna, in this relation and to thank her.

The book ``\emph{Plane algebroid curves in arbitrary characteristic}''~\cite{PP:22}, by Gerhard Pfister and  Arkadiusz P{\l}oski, has been a main source of inspiration for this article. Arkadiusz sent it to the second author November 2023 and it was the source of many ideas in the paper. Thanks to them both.
The second author wants to thank Michel Raibaut for interesting discussions and the two authors are very grateful to the referees for careful reading of their article.

\section{Hamburger-Noether algorithm}\label{sec:hn}

The usual Newton algorithm may not work in positive characteristic
when the weights and the characteristic are not coprime. The Hamburger-Noether algorithm works
in any characteristic, since it follows the sequence of blow-ups that solves the singularity.
The Newton algorithm can be interpreted in terms of weighted blow-ups which involve quotient
singularities by the action of some abelian group; these singularities are not well-defined when the characteristic of the field is not coprime 
with the order of the group.

\subsection{Preliminaries}
\mbox{}

Let $p,q\in\bn$ coprime with $p\geq q\geq 1$. Consider the Euclidean algorithm:
\[
p = k_1 q + r_1,\ q = k_2 r_1 + r_2,\ r_1 = k_3 r_2 + r_3,\dots\ r_{\omega-1} = k_{\omega + 1} r_{\omega} + 1,\ r_{\omega} = k_{\omega + 2}.
\]
Define $r_0:=q$, $m_1:=k_1$, $n_1:=1$, $\tilde{m}_1:=k_2$, $\tilde{n}_1:=1$, and $m_i,n_i,\tilde{m}_i,\tilde{n}_i$, $i\geq 1$ satisfying
\[
p = m_i r_{i-1} + n_i r_i,\quad q = \tilde{m}_i r_i + \tilde{n}_i r_{i+1}.
\]

\begin{lemma}
For $i\geq 1$, we have
\begin{align*}
m_{i+1} &= m_i k_{i+1} + n_i, & n_{i+1} &= m_i,\\
\tilde{m}_{i+1} &= \tilde{m}_i k_{i+2} + \tilde{n}_i, & \tilde{n}_{i+1} &= \tilde{m}_i.\\
\end{align*}
\end{lemma}

\begin{proof}
By induction.
\end{proof}

\begin{lemma}
For $i\geq 1$, $\Delta_i := n_{i + 1} \tilde{m}_i - \tilde{n}_i m_{i +1}$
equals $(-1)^i$.

\end{lemma}

\begin{proof}
We can see that $\Delta_1=-1$ and $\Delta_i+\Delta_{i-1}=0$ if $i\geq 2$.
\end{proof}

\subsection{Algorithm}
\mbox{}

Let
\begin{align*}
x &= x_1 y_1^{k_{1}}, & y &= y_1,\\
x_1 &= x_2, & y_1 &= x_2^{k_{2}} y_2,\\
x_2 &= x_3 y_3^{k_{3}}, & y_2 &= y_3.
\end{align*}
More generally,
\begin{align*}
x_i &= x_{i+1}, & y_i &= x_{i+1}^{k_{i+1}} y_{i+1},\quad i\text{ odd},\\
x_i &= x_{i+1} y_{i+1}^{k_{i+1}},& y_i &= y_{i+1},\quad i\text{ even}.\\
\end{align*}

The following lemma is proved by induction.

\begin{lemma}
For $i\geq 0$,
\begin{align*}
x &= x_{2i+1}^{n_{2i+1}}, y_{2i+1}^{m_{2i+1}},\\
y &= x_{2i+1}^{\tilde{n}_{2i}} y_{2i+1}^{\tilde{m}_{2i}};
\end{align*}
for $i\geq 1$
\begin{align*}
x &= x_{2i}^{m_{2i}} y_{2i}^{n_{2i}},\\
y &= x_{2i}^{\tilde{m}_{2i-1}} y_{2i}^{\tilde{n}_{2i-1}}.
\end{align*}
\end{lemma}

%\begin{proof}
%It can be proved by induction.
%\end{proof}

Since $r_{\omega+1}=1$, $r_{\omega+2}=0$, we have $p=m_{\omega+2}$, $q=\tilde{m}_{\omega+1}$.

\begin{lemma}
If $\omega$ is odd
\begin{align*}
x &= x_{\omega+2}^{n_{\omega+2}} y_{\omega+2}^{p},\\
y &= x_{\omega+2}^{\tilde{n}_{\omega+1}} y_{\omega+2}^{q};
\end{align*}
if $\omega$ is even
\begin{align*}
x &= x_{\omega+2}^{p} y_{\omega+2}^{n_{\omega+2}}, \\
y &=  x_{\omega+2}^{q} y_{\omega+2}^{\tilde{n}_{\omega+1}}.
\end{align*}
\end{lemma}

\begin{proof}
If $n$ is odd, $(x_{\omega+2},y_{\omega+2})=(Y_1,X_1)$,
and if $n$ is even, $(x_{\omega+2},y_{\omega+2})=(X_1,Y_1)$.
Then $x=X_1^p Y_1^{q'}$, $y=X_1^q Y_1^{p'}$,
\[
\det
\begin{pmatrix}
p& q'\\
q& p'\\
\end{pmatrix}
=1.
\qedhere
\]
\end{proof}
We can check that $p'=n_{\omega+2}\leq p$ and $q'=\tilde{n}_{\omega+1}\leq q$.
Let
\begin{equation}\label{eq:P}
P(x,y)=x^a y^b
\prod_{i=1}^k\left(
x^q - \mu_iy^p
\right)^{\nu_i},\qquad
N:= a p+ b q + p q\sum_{i=1}^k \nu_i.
\end{equation}
The sequence of Hamburger-Noether maps gives if $\omega$ is odd:
\begin{align*}
P(x,y)&:=x_{\omega+2}^{a n_{\omega+2} + b \tilde{n}_{\omega+1}} y_{\omega+2}^N
\prod_{i=1}^k\left(
x_{\omega+2}^{n_{\omega+2}q} - \mu_ix_{\omega+2}^{\tilde{n}_{\omega+1}p}
\right)^{\nu_i}\\
&=
x_{\omega+2}^\ell y_{\omega+2}^N
\prod_{i=1}^k\left(
x_{\omega+2} - \mu_i
\right)^{\nu_i},
\end{align*}
where $\ell=a n_{\omega+2} + b \tilde{n}_{\omega+1}+\tilde{n}_{\omega+1}p\sum_{i=1}^k\nu_i$.
If it is even, we obtain
\[
x_{\omega+2}^N y_{\omega+2}^\ell
\prod_{i=1}^k\left(1-\mu_i y_{\omega+2}
\right)^{\nu_i}.
\]
Let us assume that $f\in\bk[[x,y]]$ has a $(p,q)$-edge in its Newton polygon and 
its \emph{face polynomial} (see page~\pageref{facepol}) for this edge is~$P$.
Hence, if we put
\begin{equation}\label{eq:hn}
x=(Y_1+\bar{\mu})^{q'} X_1^p, \quad y=(Y_1+\bar{\mu})^{p'} X_1^q,\quad \bar{\mu}=\mu_i^{\pm 1},
\end{equation}
then
\[
f(x,y)=X_1^N \underbrace{(Y_1^{\nu_i}+\dots)}_{f_1(X_1,Y_1)}.
\]
The germ of plane curve defined by  $f_1(X_1,Y_1)$ is called
an HN-transform of~$f$.

\begin{definition}
A germ of plane curve is \emph{non-degenerate} if no HN-transfor\-mation is needed,
i.e. the exponents $\nu_i$ are all equal to~$1$.
%These Newton transformations are birational maps and their composition
%is then a birational map. Given a germ $C=\Spec\bk[[x, y]]/(f)$ for 
%$f\in\bk[[x, y]]$ its strict transform is obtained as the closure
%of the preimage of the generic point of $C$ in the source of the transformation.
%The total transform of $C$ is the divisor in the source defined as the sum 
%of the total transform and the exceptional divisors counted with multiplicities. 
\end{definition}

 \section{Newton trees}\label{sec:nt}

Let $f\in\bk[[x, y]]$ be a nonconstant reduced series, $\bk$ algebraically closed. Let us consider its Newton polygon with
$k$ edges $S_1,\dots, S_k$, ordered from top to bottom. The edge $v_i$ representing $S_i$ is supported by the line $p_i X + q_i Y = N_i$,
with $\gcd(p_i,q_i)=1$.

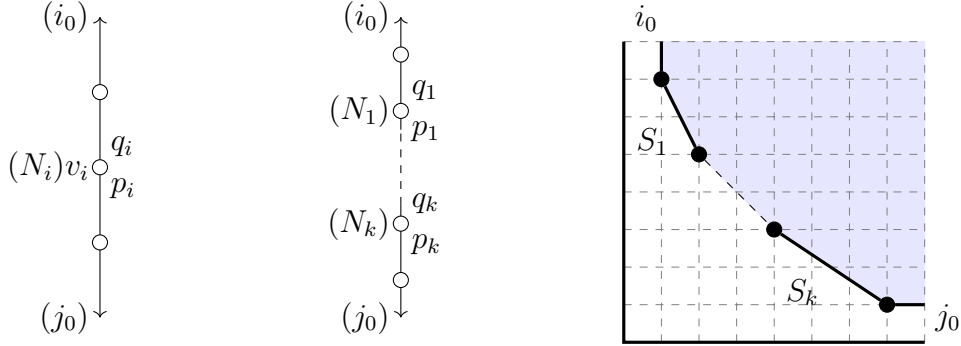
\begin{figure}[ht]
﻿\begin{tikzpicture}
\draw[<->] (0,-2) -- (0,2);
\foreach \x in {-1,0,1}
{\filldraw[fill=white] (0,\x) circle [radius=.1cm];}
\node[below right] at (0,0) {$p_i$};
\node[above right] at (0,0) {$q_i$};
\node[left] at (0,0) {$(N_i) v_i$};
\node[left] at (0,-2) {$(j_0)$};
\node[left] at (0,2) {$(i_0)$};
\begin{scope}[xshift=4cm]
\draw[<-] (0,-2) -- (0, -.5);
\draw[->](0, .5) -- (0,2);
\draw[dashed] (0, -.5) -- (0, .5);
\foreach \x in {-1.5, -.75,.75,1.5}
{\filldraw[fill=white] (0,\x) circle [radius=.1cm];}
\node[below right] at (0,-.75) {$p_k$};
\node[above right] at (0,-.75) {$q_k$};
\node[left] at (0,-.75) {$(N_k)$};
\node[below right] at (0,.75) {$p_1$};
\node[above right] at (0,.75) {$q_1$};
\node[left] at (0,.75) {$(N_1)$};
\node[left] at (0,-2) {$(j_0)$};
\node[left] at (0,2) {$(i_0)$};
\end{scope}
\end{tikzpicture}
 \hspace{2cm}
\begin{tikzpicture}[scale=.5,vertice/.style={draw,circle,fill,minimum size=0.2cm,inner sep=0}]
\draw[fill, color=blue!10!white] (1,8) --(1,7)--(2,5)--(4,3)--(7,1)--(8,1)--(8,8)--cycle;
\draw[help lines,step=1,dashed] (0,0) grid (8,8);
\draw[very thick] (0,8) -- (0,0) --( 8,0);

\node[vertice] (a17) at (1,7) {};
\node[vertice] (a25) at (2,5) {};
\node[vertice] (a43) at (4,3) {};
\node[vertice] (a70) at (7,1) {};

\draw[,very thick] ($(a17)+(0,1)$) --(a17) -- node[below left] {$S_1$} (a25)  (a43)--node[below left] {$S_k$} (a70)--($(a70)+(1,0)$);
\draw[dashed] (a25) -- (a43);
\node[above right] at (0,8) {$i_0$};
\node[above right] at (8,0) {$j_0$};
\end{tikzpicture}
 \caption{\footnotesize Left figure corresponds to the part of the Newton tree associated to the edge $v_i$. The central one is the Newton tree of the right figure.}
\label{fig:vi}
\end{figure}
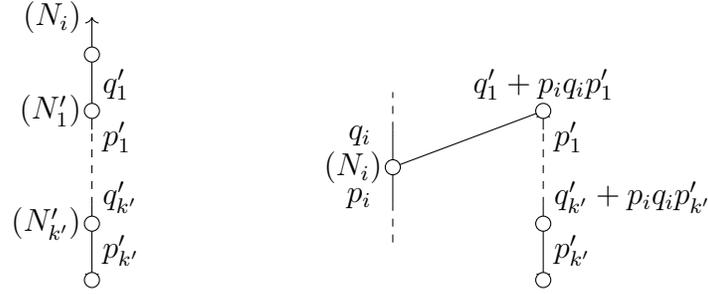
\begin{figure}[ht]
﻿\begin{tikzpicture}
\begin{scope}[xshift=0cm]
\draw[<-] (0,-1.5) -- (0, -.5);
\draw[->](0, .5) -- (0,2);
\draw[dashed] (0, -.5) -- (0, .5);
\foreach \x in {-1.5, -.75,.75,1.5}
{\filldraw[fill=white] (0,\x) circle [radius=.1cm];}
\node[below right] at (0,-.75) {$p'_{k'}$};
\node[above right] at (0,-.75) {$q'_{k'}$};
\node[left] at (0,-.75) {$(N'_{k'})$};
\node[below right] at (0,.75) {$p'_{1}$};
\node[above right] at (0,.75) {$q'_{1}$};
\node[left] at (0,.75) {$(N'_1)$};
\node[left] at (0,2) {$(N_i)$};
\end{scope}
\begin{scope}[xshift=4cm]
\draw (0,0) -- (2,.75);
\draw (0,-.5) -- (0,.5);
\draw[dashed] (0, -1) -- (0, -.5) (0, .5) -- (0, 1);
\foreach \x in {0}
{\filldraw[fill=white] (0,\x) circle [radius=.1cm];}
\node[below left=4pt] at (0,0) {$p_i$};
\node[above left=4pt] at (0,0) {$q_i$};
\node[left] at (0,0) {$(N_i)$};

\begin{scope}[xshift=2cm]
\draw[<-] (0,-1.5) -- (0, -.5);
\draw[->](0, .5) -- (0,.75);
\draw[dashed] (0, -.5) -- (0, .5);
\foreach \x in {-1.5, -.75,.75}
{\filldraw[fill=white] (0,\x) circle [radius=.1cm];}
\node[below right] at (0,-.75) {$p'_{k'}$};
\node[above right] at (0,-.75) {$q'_{k'}+p_i q_i p'_{k'}$};

\node[below right] at (0,.75) {$p'_{1}$};

\node[above] at (0,.75) {$q'_1+p_i q_i p'_1$};

\end{scope}
\end{scope}
\end{tikzpicture}
 \caption{\footnotesize New Newton tree and gluing.}
\label{fig:pegado}
\end{figure}

The Newton polygon is represented by a vertical linear tree. 
%and the decorations satisfy
%\[
%\frac{p_1}{q_1}>\dots>\frac{p_k}{q_k}.
%\]
Each face is represented by a vertex of the tree, and two
vertices of the tree are connected by an edge if and only if the corresponding faces intersect.
Each edge is decorated at its extremities by natural numbers. Near a vertex $v_i$, representing
a face with equation $p_i X + q_i Y = N_i$, the edges arising from $v_i$
are decorated with $p_i$ and $q_i$ (see Figure~\ref{fig:vi}). The vertex $v_i$ is decorated
with $N_i$.
If $x^{i_0}$ and $y^{j_0}$ are factors  of $f$ (with maximal multiplicity, $0\leq i_0, j_0\leq 1$) then the
non-compact faces $X=i_0, Y=j_0$ are represented by two decorated arrows.
%, see Figure~\ref{fig:vi}, where we can see the decorations near the vertices
%reflecting $p_i,q_i,N_i$ (note the positions of the two first ones). 
The decorated tree 
contains the same information as the Newton polygon.
As we did in \eqref{eq:P}, each face has an associated homogeneous polynomial (the \emph{face polynomial})\label{facepol}
\begin{equation*}
P_i(x,y)=x^{n_i} y^{m_i}
\prod_{j=1}^{k_i}\left(
x^{q_i} - \mu_{i,j}y^{p_i}
\right)^{\nu_{i,j}},\quad
\gcd(p_i, q_i)=1.
\end{equation*}
Let us fix $i,j$; if $\nu_{i,j}=1$, we attach to the vertex $v_i$ an edge ended with an arrow (to the right). If $\nu_{i,j}>1$,
we apply Hamburger-Noether algorithm \eqref{eq:hn} for $(p_i,q_i,\mu_{i,j})$ (the choice of the power $\pm1$ depends on the
parity of the length of the Euclidean algorithm):
\[
f_{i,j}(x_1, y_1) = x_1^{N_i}(y_1^{\nu_{i,j}}+\dots)\in\bk[[x,y]]
\]
to which a new Newton polytope can be applied (in general with a smaller height) and translated to the right by $N_i$. We glue the new Newton tree
as in Figure~\ref{fig:pegado}, including the changes of decorations.

\begin{remark}
The construction of Newton trees in characteristic~$0$ is similar, but we can use Newton maps instead of Hamburger-Noether maps, with the same result. 
\end{remark}

\subsection{Examples}
\mbox{}

\begin{example}\label{ex1}
Let
\[
f(x,y) := (x^2 - y^3)^4 -2 (x^2 - y^3)^2 x y^{11}- y^{19} (1 - y^3)  (x^2 - y^3) + y^{25}.
\]
The Newton polygon has only one edge with face polynomial $P(x,y):=(x^2 - y^3)^4$
with $\omega=0$.
There is only one edge and one root, hence only one transformation:
\[
x=x_1^3 (y_1 + 1),\qquad
y=x_1^2 (y_1 + 1).
\]
Then
\[
f(x, y) = x_1^{24} (y_1 + 1)^8 ((y_1^2 - x_1^{13})^2 + \dots).
\]
The Newton polytope is associated to $(y_1^2 - x_1^{13})^2$, $\omega=0$.
Hence,
\[
x_1=x_2^2 (y_2 + 1),\qquad
y_1=x_2^{13} (y_2 + 1)^6,
\]
and
\[
f(x, y) = x_2^{100} (y_2 + 1)^{48} (1 + x_2^{13} + \dots)^8 (y_2^2 + x_2 + \dots).
\]
% The next Newton polytope is associated to $5776 y_2^2 - x_2$, $\tilde{n}=1$.
% Hence,
% \[
% x_2=x_3^2 (y_3 + 5776),\qquad
% y_2=x_3.
% \]
% and
% \[
% f(x, y) = x_3^{202} u(x_3, y_3) (y_3 - 242592 x_3 + \dots), \qquad u(0,0)\neq 0.
% \]
\begin{figure}[ht]
\begin{tikzpicture}
\foreach\x in {1,...,4}
{
\coordinate (A\x) at (2*\x, 0);
}
\draw[->] (A1) -- (A4);
\draw[->] (A1) -- ($(A1) + (0,1)$) node[above] {$(0)$};
\foreach\x in {1,...,3}
{
\draw[->] (A\x) -- ($(A\x) - (0,1)$) node[below] {$(0)$};
}
\foreach\x in {1,...,3}
{
\filldraw[fill=white] (A\x) circle [radius=.1];
}
\node[above left] at (A1) {$3$};
\node[below left] at (A1) {$2$};
\node[left=3pt] at (A1) {$(24)$};
\node[above left] at (A2) {$25$};
\node[below right] at (A2) {$2$};
\node[above=10pt] at (A2) {$(100)$};
\node[above left] at (A3) {$101$};
\node[below right] at (A3) {$2$};
\node[above=10pt] at (A3) {$(202)$};
\end{tikzpicture}
 \caption{\footnotesize Newton tree of Example~\ref{ex1}.}
\label{fig:ej1}
\end{figure}
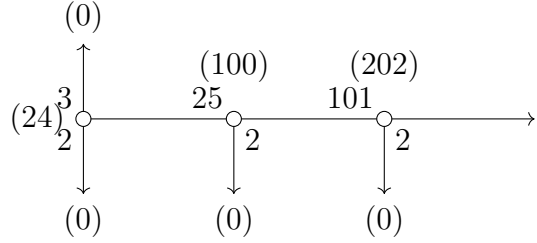

\end{example}

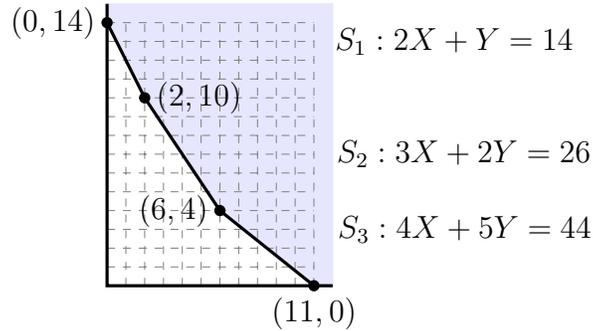
\begin{figure}[ht]
﻿\begin{tikzpicture}[scale=.25,vertice/.style={draw,circle,fill,minimum size=0.2cm,inner sep=0}]
\foreach \x/\y in {0/14, 2/10, 6/4, 11/0}
{
\coordinate (A\x-\y) at (\x,\y);
}

\draw[fill, color=blue!10!white] ($(A0-14)+(0,1)$) --(A0-14)--(A2-10)--(A6-4)--(A11-0)--($(A11-0)+(1,0)$)--(12,15)--cycle;
\draw[help lines,step=1,dashed] (0,0) grid (11,14);
\draw[very thick] (0,15) -- (0,0) --(12,0);

\fill (A0-14) circle [radius=.3];
\fill (A2-10) circle [radius=.3];
\fill (A6-4) circle [radius=.3];
\fill (A11-0) circle [radius=.3];

\draw[,very thick] (A0-14) node[left] {$(0,14)$} --
node[right=2.75cm, pos=.25] {$S_1: 2X+Y=14$}
(A2-10) node[right] {$(2,10)$}
--
node[right=1.9cm, pos=.5] {$S_2: 3X+2Y=26$}
(A6-4)  node[left] {$(6,4)$} --
node[right=1.1cm, pos=.25] {$S_3: 4X+5Y=44$}
(A11-0)  node[below] {$(11, 0)$};

\end{tikzpicture}
 \caption{\footnotesize Newton polygon of Example~\ref{ex2}.}
\label{fig:ej2a}
\end{figure}

\begin{example}\label{ex2}
Let
\[
f(x,y) := -x^2 y^4 (x^2 - y^3)^2 +x^{11} + y^{14}  + x y^{13}.
\]
The factorized face polynomials are
\[
P_1=
\begin{cases}
y^{10} (y^2 -x)(y^2+x)&\text{ if }\chr\bk\neq 2,\\
y^{10} (y^2 -x)^2&\text{ if }\chr\bk= 2;
\end{cases}
\]
$P_2=x^2 y^4(x^2 - y^3)^2$,
and
$P_3=x^6(x^5 - y^4)$.
For~$S_1$ and $\chr\bk=2$, we take 
\[
x=x_1^2 (y_1 + 1),\qquad
y=x_1,
\]
for which we obtain that, up to a unit in $\bk[[x_1,y_1]]$, $f$ is $x_1^{14}  (y_1^2 + x_1 +\dots)$.
For~$S_2$:
\[
x=x_1^3 (y_1 + 1),\qquad
y=x_1^2 (y_1 + 1).
\]
Then
\[
f(x, y) = x_1^{26} (y_1 + 1)^{10} (x_1^2 - y_1^{2} + \dots).
\]
If $\chr\bk=2$, we perform the transformation
\[
x_1=x_2 (y_2 + 1),\qquad
y_1=x_2,
\]
and we obtain for $f$, up to a unit in $\bk[[x_2,y_2]]$, $x_2^{28}  (y_2^2 + x_2 +\dots)$.
If $\chr\bk\neq2$, we finish with two arrows.

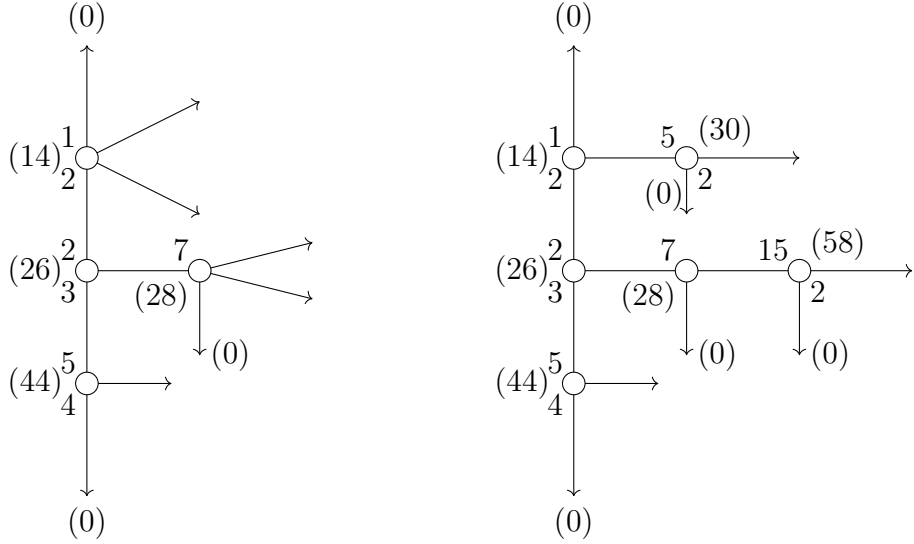
\begin{figure}[ht]
\begin{tikzpicture}[scale=1.5]
\foreach\x in {0,1,...,4}
{
\coordinate (A\x) at (0,\x);
}
\draw[<->] (A0) node[below] {$(0)$} -- (A4)node[above] {$(0)$};
\draw[->] (A1) -- ($(A1) + (.75,0)$);
\node[above left] at (A1) {$5$};
\node[below left] at (A1) {$4$};
\node[left=5pt] at (A1) {$(44)$};

\draw (A2) -- ($(A2) + (1,0)$) node[above left] {$7$}  node[below left] {$(28)$};
\draw[->] ($(A2) + (1,0)$) -- ($(A2) + (1,-.75)$)  node[right] {$(0)$};
\node[above left] at (A2) {$2$};
\node[below left] at (A2) {$3$};
\node[left=5pt] at (A2) {$(26)$};
\draw[->] ($(A2) + (1,0)$) -- ($(A2) + (2,.25)$);
\draw[->] ($(A2) + (1,0)$) -- ($(A2) + (2,-.25)$);

\draw[->] ($(A3)$) -- ($(A3) + (1,.5)$);
\draw[->] ($(A3)$) -- ($(A3) + (1,-.5)$);
\node[above left] at (A3) {$1$};
\node[below left] at (A3) {$2$};
\node[left=5pt] at (A3) {$(14)$};

\foreach\x in {1,...,3}
{
\filldraw[fill=white] (A\x) circle [radius=.1];
}
\filldraw[fill=white] ($(A2) + (1,0)$) circle [radius=.1];

\end{tikzpicture}
 \hspace{2cm}
\begin{tikzpicture}[scale=1.5]
\foreach\x in {0,1,...,4}
{
\coordinate (A\x) at (0,\x);
}
\draw[<->] (A0) node[below] {$(0)$} -- (A4)node[above] {$(0)$};
\draw[->] (A1) -- ($(A1) + (.75,0)$);
\node[above left] at (A1) {$5$};
\node[below left] at (A1) {$4$};
\node[left=5pt] at (A1) {$(44)$};

\draw (A2) -- ($(A2) + (2,0)$) ;
\node[above left] at (A2) {$2$};
\node[below left] at (A2) {$3$};
\node[left=5pt] at (A2) {$(26)$};
\draw[->] ($(A2) + (1,0)$) node[above left] {$7$} node[below left] {$(28)$} -- ($(A2) + (1,-.75)$) node[right] {$(0)$};
\draw[->] ($(A2) + (2,0)$) -- ($(A2) + (3,0)$);
\draw[->] ($(A2) + (2,0)$) node[above left] {$15$} node[above right] {$(58)$} node[below right] {$2$} -- ($(A2) + (2,-.75)$) node[right] {$(0)$};

\draw[->] ($(A3)$) -- ($(A3) + (2,0)$);
\draw[->] ($(A3)+(1,0)$) node[above left] {$5$}
node[below right] {$2$} node[above right] {$(30)$}
-- ($(A3) + (1,-.5)$) node[above left=-3pt] {$(0)$};
\node[above left] at (A3) {$1$};
\node[below left] at (A3) {$2$};
\node[left=5pt] at (A3) {$(14)$};

\foreach\x in {1,...,3}
{
\filldraw[fill=white] (A\x) circle [radius=.1];
}
\filldraw[fill=white] ($(A2) + (1,0)$) circle [radius=.1];
\filldraw[fill=white] ($(A2) + (2,0)$) circle [radius=.1];
\filldraw[fill=white] ($(A3) + (1,0)$) circle [radius=.1];

\end{tikzpicture}
 \caption{\footnotesize Newton tree of Example~\ref{ex2} if $\chr\bk\neq 2$ (to the left),
or $\chr\bk= 2$ (to the right).}
\label{fig:ej2}
\end{figure}
\end{example}

 \subsection{Minimal trees}
\mbox{}

Let us consider a Newton tree $\cT$. Let $\cV$ be its set of vertices, $\cE$ its set of edges, $\arista$ its set of arrows,
and $\arista_0$ the set of its arrows decorated by~$(0)$. 
\begin{definition}
A \emph{dead end} is an edge between a vertex and an arrow decorated by $(0)$.
\end{definition}

We will proceed with the following conventions:

\begin{enumerate}[label=(\alph{enumi}), leftmargin=*, widest=b]
\item\label{mov-a} The dead ends decorated by $1$ and the attached arrows
will be erased.
\item\label{mov-b} Vertices of valency~$2$  will be erased, keeping the decorations of the remaining vertices.
\end{enumerate}

\begin{figure}[ht]
\centering
\begin{tikzpicture}
\begin{scope}
\draw[dashed] (-1, 0) -- (-.5,0) (.5,0) -- (1,0);
\draw (-.5,0) -- (.5, 0);
\draw[->] (0,0) -- (0,-.75) node[below] {$(0)$};
\filldraw[fill=white] (0,0) node[above=2pt] {$v$} circle [radius=.1cm];
\end{scope}
\begin{scope}[xshift=3cm]
\draw[dashed] (-1, 0) -- (-.5,0) (.5,0) -- (1,0);
\draw (-.5,0) -- (.5, 0);
\filldraw[fill=white] (0,0) node[above=2pt] {$v$}  circle [radius=.1cm];
\end{scope}
\begin{scope}[yshift=-2cm, xshift=-1cm]
\draw[dashed] (-2, 1) -- (-1.5,.5) (-2, -1) -- (-1.5,-.5);
\draw (-1.5,.5) -- (-1, 0) -- (-1.5,-.5) (-1, 0) -- (-1,-.5);
\draw (-1,0) -- (1, 0);
\draw[xscale=-1, dashed] (-2, 1) -- (-1.5,.5) (-2, -1) -- (-1.5,-.5);
\draw[xscale=-1] (-1.5,.5) -- (-1, 0) -- (-1.5,-.5) (-1, 0) -- (-1,-.5);
\draw (-1,0) -- (1, 0);
\draw[->] (0,0) -- (0,-.75) node[below] {$(0)$};
\filldraw[fill=white] (0,0) node[above=2pt] {$v$}  circle [radius=.1cm];
\filldraw[fill=white] (-1,0) node[above=2pt] {$v'$}  circle [radius=.1cm];
\filldraw[fill=white] (1,0) node[above=2pt] {$v''$}  circle [radius=.1cm];
\end{scope}
\begin{scope}[yshift=-2cm, xshift=4cm]
\draw[dashed] (-2, 1) -- (-1.5,.5) (-2, -1) -- (-1.5,-.5);
\draw (-1.5,.5) -- (-1, 0) -- (-1.5,-.5) (-1, 0) -- (-1,-.5);
\draw (-1,0) -- (1, 0);
\draw[xscale=-1, dashed] (-2, 1) -- (-1.5,.5) (-2, -1) -- (-1.5,-.5);
\draw[xscale=-1] (-1.5,.5) -- (-1, 0) -- (-1.5,-.5) (-1, 0) -- (-1,-.5);
\draw (-1,0) -- (1, 0);

\filldraw[fill=white] (-1,0) node[above=2pt] {$v'$}  circle [radius=.1cm];
\filldraw[fill=white] (1,0) node[above=2pt] {$v''$}  circle [radius=.1cm];
\end{scope}
\end{tikzpicture}
 \caption{\footnotesize Erasing operations.}
\label{fig:borrado}
\end{figure}
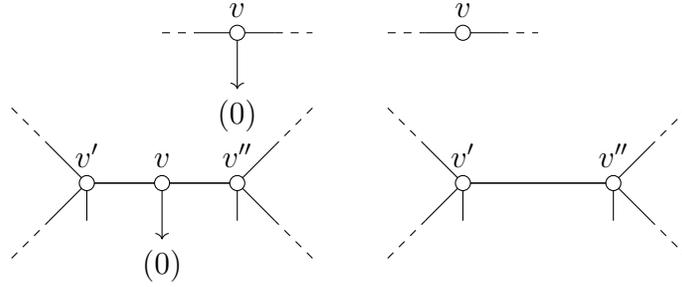

\begin{remark} These operations are defined in \cite[Def.~1.1.1]{cnd}, and their properties are studied. A tree is minimal if no operation \ref{mov-a}-\ref{mov-b} can be performed.
\end{remark}

\begin{example}
The tree in Figure~\ref{fig:ej1} is already minimal. In Figure~\ref{fig:ej2}, left, the dead end to the left can be erased and the tree becomes minimal.
The minimalization of Figure~\ref{fig:ej2}, right, is in Figure~\ref{fig:min}.
\begin{figure}[ht]
\centering
﻿\begin{tikzpicture}[scale=1.5]
\foreach\x in {0,1,...,4}
{
\coordinate (A\x) at (0,\x);
}
\draw[<->] (A0) node[below] {$(0)$} -- (A4)node[above] {$(0)$};
\draw[->] (A1) -- ($(A1) + (.75,0)$);
\node[above left] at (A1) {$5$};
\node[below left] at (A1) {$4$};
\node[left=5pt] at (A1) {$(44)$};

\draw (A2) -- ($(A2) + (2,0)$) ;
\node[above left] at (A2) {$2$};
\node[below left] at (A2) {$3$};
\node[left=5pt] at (A2) {$(26)$};

\draw[->] ($(A2) + (1,0)$) -- ($(A2) + (2,0)$);
\draw[->] ($(A2) + (1,0)$) node[above left] {$15$}  node[below right] {$2$} node[above right] {$(58)$} -- ($(A2) + (1,-.75)$) node[right] {$(0)$};

\draw[->] ($(A3)$) -- ($(A3) + (1,0)$);

\node[above left] at (A3) {$2$};
\node[below left] at (A3) {$5$};
\node[left=5pt] at (A3) {$(30)$};

\foreach\x in {1,...,3}
{
\filldraw[fill=white] (A\x) circle [radius=.1];
}
\filldraw[fill=white] ($(A2) + (1,0)$) circle [radius=.1];

\end{tikzpicture}
 \caption{\footnotesize}
\label{fig:min}
\end{figure}
\end{example}

\begin{definition}
Let $v\in\cV$. The \emph{valency} $\delta_v$ of $v$ is the number of edges through~$v$ (the valency of an arrow is~$1$). Let $v'\in\arista_0$ be attached to $v\in\cV$ such that $a$ is the decoration of the
edge joining $v,v'$. Then, the \emph{multiplicity} of $v'$ is $N_{v'}:=\frac{N_v}{a}$. Finally, the \emph{multiplicity} of $\cT$ is
\[
M(\cT):= -\sum_{v\in\cV\cup\arista_0} N_v(\delta_v-2),
\]
where $\delta_v$ is the valency of $v$ in the tree.
\end{definition}

\begin{example}
For Example~\ref{ex1}, the multiplicity is $M(\cT)=-155$; for Example~\ref{ex2}, it is $-103$ if $\chr\bk= 2$, and $-101$ otherwise.
\end{example}

 \section{Multiplicity of a tree and area of Newton polygons}\label{sec:area}

\begin{definition}
A \emph{convenient Newton polygon} is a Newton polygon which hits both axes.Let $\area$ be the area between the Newton polygon and the axes. 

A semi-convenient Newton polygon is a Newton polygon that hits the line $y=0$ and the line $x=N$.  In this case, let $\area$ be the area between the Newton polygon and the lines $x=N$ and $y=0$. 

\end{definition}

\begin{definition}
A \emph{reduced convenient Newton polygon} is a Newton polygon which hits the lines $x=0$ or $x=1$ and $y=0$ or $y=1$.  
\end{definition}

\begin{remark}
We can make such Newton polygon convenient in case it does not hit $x=0$ or $y=0$. This is done as follows. Assume that the polygon does not hit $x=0$.
Let $v$ be the vertex of the Newton polygon $\mathcal{N}$ with coordinates $(1,\beta)$ and 
let $v_n$ be a point with coordinates $(0,n) $ with $n$ large enough such that the set of vertices of the  convex hull of $\cV':=\cV\cup\{v_n\}$  is~$\cV'$.
A similar procedure can be applied if the polygon
does not hit $y=0$.  If $a$ and $b$ are the lengths on the axes, $2\area-a-b$ does not depend on how we make the Newton polygon convenient.
\end{remark}

\begin{definition}
A \emph{reduced semi-convenient Newton polygon} is a Newton polygon which hits the lines $x=N $ and $y=0$ or $y=1$. Again, we can make it semi-convenient and if $a+N$ is the length on the $x$-axis, $2\area-a$ does not depend on the chosen semi-convenient Newton polygon.
\end{definition}

\begin{definition}
A tree is \emph{non-degenerate} if it is \emph{vertical} except for the non-decorated arrows. For the associated Newton polygon, the face polynomials are as in \eqref{eq:P} with $\nu_i=1$.
In particular, the Hamburger-Noether algorithm is not needed.
\end{definition}

\begin{remark}
In positive characteristic this definition does not imply that the face polynomials do not 
have singularities in the torus. 
%Both definitions do not coincide
%in positive characteristic. 
For example, let
\[
f(x,y):= x y (x + y),
\]
which is non-degenerate according to this definition but $(1,1)$ is a common zero of the derivatives in $\mathbb{F}_3^2$.
\end{remark}

\begin{lemma}[{\cite[Lemma 5.20]{cnv:14}}]
If the tree $\cT$ is reduced, and non-degenerate, then
\[
-M(\cT)=2 \area - a - b.
\]
\end{lemma}

\begin{theorem}\label{thm:cv14}
Let $f\in\bk[[x,y]]$, and let $\cT$ be its Newton tree. Then,
\[
-M(\cT)=2 \area_0 - a - b + \sum_{\ell=1}^r  (2\area_\ell-a_\ell),
\]
where $\area_0$ is the area of the first Newton polygon and $a,b$
its traces on the axes. The summation is taken over all Hamburger-Noether
transforms, $\area_\ell$ is the area of the Newton polygon limited by the line $x=N_\ell$
and $a_\ell+N_\ell$ is the trace on the $x$-axis.
\end{theorem}

\begin{proof}
 Let $v$ be a vertex, and let  $AB$ be the face of the Newton polygon which corresponds to $v$ in the HN algorithm. Let  $\delta_v$ be the valency of $v$, and let $d_v$ be the number of points with integral coordinates on the face $AB$.
 
 \begin{figure}[ht]
 \centering
 \begin{tikzpicture}
\draw (0, 2) -- (0,0) node[below left] {$O$} -- (2,0);
\draw (0,0) -- (.5,1.5) node [above] {$B$} -- (1.5, .5)  node [right] {$A$} -- cycle;
\fill (.5,1.5) circle [radius=.1cm];
\fill (1.5, .5) circle [radius=.1cm];
\end{tikzpicture}
  \caption{\footnotesize}
 \label{fig:triangle}
 \end{figure}

 Let $\area_v$ be the area of the triangle $OAB$. We have
 $2\area_v-a-b=N_v(d_v-2)$
 where $a$ is the first coordinate of $A$ if $A$ is on the $x$-axis and $0$ otherwise, and $b$ the second coordinate of $B$ if $B$ is on the y-axis, and $0$ otherwise. Since
 \[
 d_v=\sum_{i=1}^{\delta_v} \nu_i 
 \Longrightarrow  d_v-\delta_v = \sum_{i=1}^{\delta_v} (\nu_i-1),
 \]
 then
 \[
 N_v(\delta_v-1) = N_v(d_v-1) - N_v(d_v-\delta_v)= N_v(d_v-1) - N_v \sum_{i=1}^{\delta_v} (\nu_i-1).
 \]
 If $\nu_i=1$ for all $i$, then  $2A_v-a-b=N_v(d_v-2)=N_v(\delta_v-2)$.
 
  \begin{figure}[ht]
 \centering
 \begin{subfigure}[b]{.49\textwidth}
  \centering
   \begin{tikzpicture}
\draw (0, 2) -- (0,0) -- (5,0);
\draw (2, 2) -- (2,0) node[above left] {$N_{v}$};
\draw (2,1.5) node[left] {$\nu_i$} -- (3,.5) -- (4.5,0);
\draw[<->] (2,-.25) -- node[below] {$a_{v,i}$} (4.5,-.25);
\draw[->] (3.25,1) node[above] {$\area_{v,i}$} -- (2.5,.5);
\draw[pattern=north east lines] (2,1.5) -- (2,0) -- (4.5,0) -- (3, .5) --cycle ;
\end{tikzpicture}
  \caption{\footnotesize }
 \label{subfig:shift}
 \end{subfigure}
  \begin{subfigure}[b]{.49\textwidth}
  \centering
 \begin{tikzpicture}
\draw (0, 2) -- (0,0) node[below left] {$O$} -- (5,0);
\draw (2, 1.5) -- (2,0) node[below] {$N_{v}$};
\draw (0,0) -- (2,1.5) node[left] {$\nu_i$} -- (3,.5) -- (4.5,0);
\draw[<->] (0,-.55) -- node[below] {$\overline{a}_{v,i}$} (4.5,-.55);
\draw[->] (3.25,1) node[above] {$\overline{\area}_{v,i}$} -- (2.5,.5);
\draw[pattern=north east lines] (2,1.5) -- (0,0) -- (4.5,0) -- (3, .5) --cycle ;
\end{tikzpicture}
  \caption{}
 \label{subfig:shift1}
  \end{subfigure}
  \caption{\footnotesize}
  \label{fig:areas}
 \end{figure}

 If there exist $i$ such that $\nu_i>1$, we consider all HN-transforms, and we compute 
 \[
 \area_{v,1}:=2 \area_v-a-b+ \sum_{i=1}^{\delta_v} (2\area_{v,i} -a_{v,i}).
 \]
Then $\area_{v,1}=N_v(\delta_v-2)+ \sum_i 2\area_{v,i}+N_v\nu_i -(N_v+a_{v,i})$. Let 
 $\overline{\area}_{v,i}:=2\area_{v,i}+N_v\nu_i$ be the area of the polygon delimited by the Newton polygon $\mathcal{N}_{v,i}$ and $O$, and let
 $\overline{a}_{v,i}:=N_v+a_{v,i}$ be its trace on the $x$-axis.

 Then we proceed by induction, and we stop when all $\nu$'s are equal to~$1$.
 \end{proof}

 \section{Newton trees and intersection multiplicity}\label{sec:inter}

In this section, we will show how to compute the intersection multiplicity of two branches using the trees.

\begin{ntc}
Let $v$ be a vertex of $\cT$, and let $\varepsilon_{v}$ be the set of edges incident to $v$. If $e \in \varepsilon_{v}$ let $q(v,e)$ be the decoration on $e$ near $v$.
Let $v,w$ be two vertices or arrows of $\cT$; we denote 
the path from $v$ to
$w$ in $\cT$ as $\varepsilon=\varepsilon_{v,w}$; $\cV_\varepsilon$ is the set of vertices and arrows in $\varepsilon$,
and $\cE_\varepsilon$ is the set of  edges in~$\varepsilon$. If $x$ is a vertex, then
$\cE_{x}$ is the set of edges containing~$x$. We denote
\begin{align*}
\rho(v,w) &:= \prod_{x\in\cV_\varepsilon\setminus\{v,w\}}
\prod_{e\in\cE_{x}\setminus\cE_\varepsilon}q(x,e),
\\
\overline{\rho}(v,w) &:= \prod_{x\in\cV_\varepsilon\setminus\{w\}}
\prod_{e\in\cE_{x}\setminus\cE_\varepsilon}q(x,e).
\end{align*}
In case $v,w$ are two distinct arrows %$\alpha,\beta$ 
we put $i(v,w):=\rho(v,w)$.  For $\cV_1,\cV_2\subset\arista$ we set
\[
i(\cV_1,\cV_2):=\sum_{v_1\in\cV_1}\sum_{v_2\in\cV_2\setminus\{v_1\}} i(v_1, v_2).
\] 
\end{ntc}

\begin{figure}[ht]
\centering
\begin{tikzpicture}
\draw (0,0) -- (1.5,0) (2.5,0) -- (3.5,0) (4.5,0) -- (6,0);
\draw[dotted] (1.5,0) -- (2.5,0) (3.5,0) -- (4.5,0);

\draw (1,0) -- (1,-.75) (5,0) -- (5,-.75);

\draw (0,0) -- (135:.75) (0,0) -- (-135:.75);
\node at (-.5,.05) {$\vdots$};

\draw[shift={(3,0)}] (0,0) --  node[right] {$e_{ij}$} (-45:.75) (0,0) -- (-135:.75);
\node at (3.05,-.5) {$\dots$};

\filldraw[fill=white] (0,0) node[above=2pt] {$v$} circle [radius=.1];

\filldraw[fill=white] (1,0) circle [radius=.1];

\filldraw[fill=white] (3,0) node[above=2pt] {$x_{i}$} circle [radius=.1];

\filldraw[fill=white] (6,0) node[above=2pt] {$w$} circle [radius=.1];

\end{tikzpicture}
 \caption{\footnotesize}
\label{fig:path}
\end{figure}

An irreducible $f\in\bk[[x,y]]$ admits primitive \emph{parametrizations}
$t\mapsto(\varphi(t),\psi(t))$. If $g$ is another series, then
\[
i(f,g):=\ord g(\varphi(t),\psi(t)).
\]
If $g$ is irreducible, we can exchange the roles of $f,g$, and the value does not change, and the definition can be
extended to the case where $f,g$ are reducible assuming that $i(f,f)=\infty$.

\begin{theorem}
Let $f,g\in\bk[[x,y]]$ be irreducible series, let $\cT$ be the Newton tree of $fg$, and let $\alpha_f,\alpha_g$
be the arrows representing the last HN-transforms of $f,g$, respectively. Then
\[
i(f,g)=i(\alpha_f,\alpha_g).
\]
\end{theorem}

\begin{proof}
We decompose the proof into several steps.

\begin{enumerate}[label=(\arabic{enumi}), leftmargin=*, widest=2]
\item Either the Newton polygons of $f,g$ are distinct or their face polynomials do not coincide.
%The arrows separate on the first Newton polygon.

\begin{enumerate}[label=(\alph{enumii}), ref=(\arabic{enumi}\alph{enumii}), leftmargin=*, widest=b]
\item\label{case1a} They separate in different vertices as in Figure~\ref{fig:inter} left.
Under these conditions, we can assume that
\[
f(x,y) = (x^{q_i} - \mu y^{p_i})^{\nu_i}+\ldots,
\]
and a parametrization of $g=0$ is given by
\[
\varphi(t)= t^{p_j \nu_j}\tilde{\varphi}(t),\quad
\psi(t)= t^{q_j \nu_j}\tilde{\psi}(t),\qquad
\tilde{\varphi}(0),\tilde{\psi}(0)\text{ units.}
\]
Then
\[
f(\varphi(t),\psi(t)) = (t^{p_j q_i \nu_j}\tilde{\varphi}^{p_i}(t)-\mu t^{q_j p_i \nu_j}\tilde{\psi}_1^{q_i}(t))^{\nu_i} + \ldots\ .
\]
Since $p_i q_j > q_i p_j$, we have that $i(f,g)=  q_i p_j\nu_i\nu_j$. On the other hand,
one can see that the product of the weights on the horizontal part of $\varepsilon$ are $\nu_i$, $\nu_j$
and the product on the vertical part is $q_i p_j$, i.e., $i(\alpha_f,\alpha_g)$.

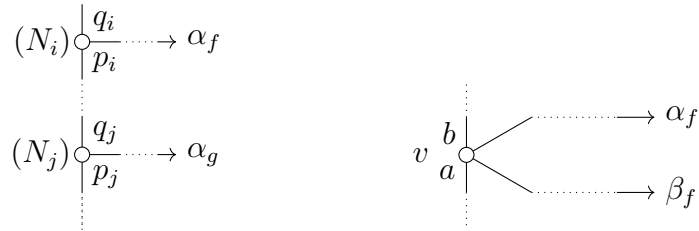
\begin{figure}[ht]
\begin{tikzpicture}

\draw[dotted] (0,-1.75) -- (0, -1.25);
\draw (0,-1.25) -- (0, -.25);
\draw(0, .25) -- (0,1.25);
\draw[dotted] (0,-1.25) -- (0, -1.75);
\draw[dotted] (0, -.25) -- (0, .25);

\draw (0,.75) -- (.5,.75);
\draw[dotted] (.5, .75) -- (1, .75);
\draw[->] (1,.75) -- (1.25, .75) node[right] {$\alpha_f$};

\draw (0,-.75) -- (.5,-.75);
\draw[dotted] (.5, -.75) -- (1, -.75);
\draw[->] (1,-.75) -- (1.25, -.75) node[right] {$\alpha_g$};

\foreach \x in {-.75,.75}
{\filldraw[fill=white] (0,\x) circle [radius=.1cm];}
\node[below right] at (0,-.75) {$p_j$};
\node[above right] at (0,-.75) {$q_j$};
\node[left] at (0,-.75) {$(N_j)$};
\node[below right] at (0,.75) {$p_i$};
\node[above right] at (0,.75) {$q_i$};
\node[left] at (0,.75) {$(N_i)$};

\end{tikzpicture}
 \hspace{2cm}
\begin{tikzpicture}

\draw (0,-.5) -- (0, .5);
\draw[dotted] (0,.5) -- (0, 1) (0,-.5) -- (0, -1);
\draw (0,0) -- (30:1) (0,0) -- (-30:1);
\draw[dotted] (30:1) -- (2, .5);
\draw[dotted] (-30:1) -- (2, -.5);
\draw[->] (2,.5) -- (2.5, .5) node[right] {$\alpha_f$};
\draw[->] (2,-.5) -- (2.5, -.5) node[right] {$\beta_f$};

\filldraw[fill=white] (0,0) node[below left] {$a$} node[above left]
{$b$} node[left=10pt] {$v$} circle [radius=.1];
\end{tikzpicture}
 \caption{\footnotesize Separations as in~\ref{case1a} and~\ref{case1b}.}
\label{fig:inter}
\end{figure}

\item\label{case1b} They separate on the same face, see~Figure~\ref{fig:inter} right. 
We have that $i=j$,  and we denote $(p_i,q_i)=:(a,b)$.
We have 
\begin{gather*}
P_{f,v}(x,y) = (x^{b} - \mu y^{a})^{\nu_1},\quad 
\varphi(t)= t^{a\nu_2}\tilde{\varphi}(t),\quad
\psi(t)= t^{q_jb\nu_2}\tilde{\psi}(t),\\
\tilde{\varphi}(0), \tilde{\psi}(0),\tilde{\varphi}(0)^b -\mu\tilde{\psi}(0)^a\neq 0,
\end{gather*}
and
\begin{gather*}
i(f,g)=\ord t^{ab\nu_1\nu_2} (\tilde{\varphi}(0)^b -\mu\tilde{\psi}(0)^{a}) =ab\nu_1\nu_2 =i(\alpha_f,\alpha_g).
\end{gather*}
\end{enumerate}

\item\label{caso2} We use induction on the number of steps such that $f, g$ separate
(see Figure~\ref{fig:caso2}). The first step has been done.
Let us assume the formula holds if we separate at step $h-1$, and we check that it holds if it separates
at step~$h$.

\begin{figure}[ht]
\centering
\begin{tikzpicture}
\draw (0,0) -- (2.5,0) (3.5, 0) -- (5,0);
\draw[dotted] (2.5,0) -- (3.5, 0);
\draw[<->] (0,-1) node[below] {$(0)$} -- (0,1) node[above] {$(0)$} ;
\draw[<->] (2,0) -- (2,-1) node[below] {$(0)$} ;
\draw[<->] (5,0) -- (5,-1) node[below] {$(0)$} ;

\draw[dotted] (6,.5) -- (6.5, .5);
\draw[->] (6.5,.5) -- (7, .5) node[right] {$\alpha_f$};
\draw[dotted] (6,-.5) -- (6.5, -.5);
\draw[->] (6.5,-.5) -- (7, -.5) node[right] {$\beta_f$};

\filldraw[fill=white] (0,0) node[above left]{$\bar{\bar{b}}$} node[below left]{$a$}  circle [radius=.1];

\filldraw[fill=white] (2,0) node[above left]{$\bar{\bar{b}}_1$} node[below right]{$a_1$} circle [radius=.1];

\filldraw[fill=white] (5,0) node[above left]{$\bar{\bar{b}}_k$} node[below right]{$a_k$}  circle [radius=.1];
\end{tikzpicture}
 \caption{\footnotesize Case~\ref{caso2}.}
\label{fig:caso2}
\end{figure}

After the first step, we have
\begin{gather*}
f(x,y)=X_1^{N_f} f_1(X_1,Y_1),\quad
g(x,y)=X_1^{N_g} g_1(X_1,Y_1),\\
x= X_1^a(Y_1-\mu)^{a'},\quad
y= X_1^b(Y_1-\mu)^{b'}.
\end{gather*}
From a primitive parametrization $(\varphi_1(t), \psi_1(t))$ of $f_1(X_1, Y_1)$ we obtain a primitive 
parametrization  $(\varphi(t), \psi(t))$ of $f(X_1, Y_1)$ using the above expression. Then
\begin{align*}
i(f,g) &= \ord  g(\varphi(t), \psi(t)) = \ord \varphi_1(t)^{N_g} g_1(\varphi_1(t), \psi_1(t))\\&
=
 \ord g_1(\varphi_1(t), \psi_1(t)) + N_g \ord \varphi_1(t) =
 i(f_1, g_1) + N_g i (X_1, f_1);
\end{align*}
by symmetry, the second term equals $N_f i (X_1, g_1)$. By induction, we have $i(f_1, g_1)=\bar{b}_k a_k d_1 d_2$. 
We claim that 
\begin{equation}\label{eq:claim}
\bar{\bar{b}}_i = a\cdot b\cdot a_1^2\cdot\ldots\cdot a_{i-1}^2 a_i + \bar{b}_i.
\end{equation}
To prove this claim, note that $\bar{\bar{b}}_1 = a b a_1 + \bar{b}_1$.
We assume \eqref{eq:claim} is true until $i$ and we prove it for $i+1$:
\begin{equation*}
\bar{\bar{b}}_{i+1} = \bar{\bar{b}}_{i}\cdot a_i\cdot a_{i+1} + b_{i+1}=
a\cdot b\cdot a_1^2\cdot\ldots\cdot a_{i-1}\cdot a_i^2\cdot a_{i+1}+ \underbrace{\bar{b}_i \cdot a_i\cdot a_{i+1} + b_{i+1}}_{\bar{b}_{i+1}},
\end{equation*}
and the claim is true. Let us consider a parametrization for $g_1$:
\[
\varphi_1(t)= t^{a_1}\tilde{\varphi}_1(t),\quad
\psi_1(t)= t^{b_1}\tilde{\psi}_1(t),\qquad
\tilde{\varphi}_1(0),\tilde{\psi}_1(0)\text{ units.}
\]
The order of $\varphi_1(t)$ equals $a_1\cdot\ldots\cdot a_k\cdot d_2$ and $N_f=a\cdot b\cdot a_1\cdot\ldots\cdot a_k\cdot  d_1$.
Since 
\[
\varphi(t)= \varphi_1(t)^a\tilde{\varphi}_2(t),\quad
\psi(t)= \psi_1(t)^b \tilde{\psi}_2(t),\qquad
\tilde{\varphi}_2(0),\tilde{\psi}_2(0)\text{ units,}
\]
and $f(\varphi(t),\psi(t))=\varphi_1(t)^{N_f}f(\varphi_1(t),\psi_1(t))$,
the computation of the order gives the statement.\qedhere

\begin{figure}[ht]
\centering
\begin{tikzpicture}
\draw (0,0) -- (2.5,0) (3.5, 0) -- (5,0);
\draw[dotted] (2.5,0) -- (3.5, 0);
\draw[<->] (0,-1) node[below] {$(0)$} -- (0,1) node[above] {$N_f+N_g$} ;
\draw[<->] (2,0) -- (2,-1) node[below] {$(0)$} ;
\draw[<->] (5,0) -- (5,-1) node[below] {$(0)$} ;

\draw[dotted] (6,.5) -- (6.5, .5);
\draw[->] (6.5,.5) -- (7, .5) node[right] {$\alpha_f$};
\draw[dotted] (6,-.5) -- (6.5, -.5);
\draw[->] (6.5,-.5) -- (7, -.5) node[right] {$\beta_f$};

\filldraw[fill=white] (0,0) node[above left]{$\bar{b}_1$} node[below left]{$a_1$}  circle [radius=.1];

\filldraw[fill=white] (2,0) node[above left]{$\bar{b}_2$} node[below right]{$a_2$} circle [radius=.1];

\filldraw[fill=white] (5,0) node[above left]{$\bar{b}_k$} node[below right]{$a_k$}  circle [radius=.1];
\end{tikzpicture}
 \caption{\footnotesize Induction step for case~\ref{caso2}.}
\label{fig:caso2ind}
\end{figure}
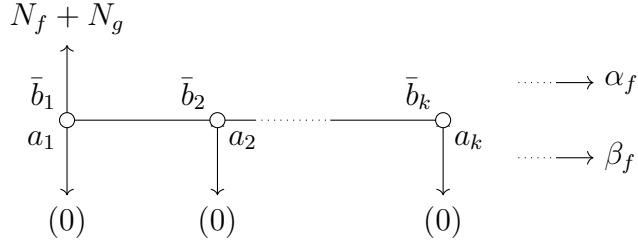

\end{enumerate}

\end{proof}
 \section{irreducible power series}\label{sec:irr}

\begin{definition}
A sequence $(v_0,v_1,\dots ,v_r)$ of positive integers is said to be a \emph{Zariski characteristic sequence} if it satisfies the following two axioms: 
\begin{enumerate}[label=\rm($Z_\arabic{enumi}$), leftmargin=*, widest=Z2]
\item Set $d_k:=\gcd (v_0,\dots ,v_k)$, $0\leq k\leq r$. Then $d_k>d_{k+1}, 0\leq k<r$, and $d_r=1$.
\item Let $n_k:=\frac {d_{k-1}}{d_k}, 1\leq k<r$. Then $n_kv_k<v_{k+1}, 1\leq k<r$.
\end{enumerate}
\end{definition}

Let $\mathcal{T}$ be a tree with one arrow not decorated by $0$. Assume that $\mathcal{T}$ is a minimal tree.
Let $v_0:=\frac{N_1}{b_1}, v_1:=\frac{N_1}{a_1},\dots ,v_r:=\frac{N_r}{a_r}$.

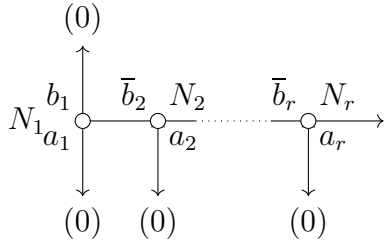
\begin{figure}[ht]
\centering
\begin{tikzpicture}
\draw[->] (0,0) -- (1.5,0) (2.5, 0) -- (4,0);
\draw[dotted] (1.5,0) -- (2.5, 0);
\draw[<->] (0,-1) node[below] {$(0)$} -- (0,1) node[above] {$(0)$} ;
\draw[<->] (1,0) -- (1,-1) node[below] {$(0)$} ;
\draw[<->] (3,0) -- (3,-1) node[below] {$(0)$} ;

\filldraw[fill=white] (0,0) node[above left]{$b_1$} node[below left]{$a_1$} node[left=10pt]{$N_1$} circle [radius=.1];

\filldraw[fill=white] (1,0) node[above left]{$\overline{b}_2$} node[below right]{$a_2$} node[above right]{$N_2$} circle [radius=.1];

\filldraw[fill=white] (3,0) node[above left]{$\overline{b}_r$} node[below right]{$a_r$} node[above right]{$N_r$} circle [radius=.1];
\end{tikzpicture}
 \caption{\footnotesize Tree $\mathcal{T}$.}
\label{fig:branche}
\end{figure}

\begin{proposition}
The sequence $(v_0,\dots ,v_r)$ is a Zariski characteristic sequence.
\end{proposition}

\begin{proof}
Note that
\begin{align*}
d_0=&v_0=\frac{N_1}{b_1}=a_1\cdot a_2\cdot\ldots\cdot a_r,\\
d_1=&\gcd(v_0,v_1)=\gcd (a_1\cdot a_2\cdot\ldots\cdot a_r,b_1\cdot a_2\cdot\ldots\cdot a_r)=a_2\cdot\ldots\cdot a_r.
\end{align*}
Then $d_1<d_0$.
Assume
\[
d_i=a_{i+1}\cdot\ldots\cdot a_r.
\]
Then
\[
d_{i+1}=\gcd (a_1\cdot a_2\cdot\ldots\cdot a_r,\dots ,\overline{b}_{i+1}\cdot a_{i+2}\cdot\ldots\cdot a_r)=a_{i+2}\cdot\ldots\cdot a_r.
\]
Then  $d_i>d_{i+1}$, and $Z_1$ is proved.
We have 
\begin{align*}
n_kv_k<v_{k+1} \iff& N_k<\frac{N_{k+1}}{a_{k+1}},\\
N_k = &\ \overline{b}_k\cdot a_k\cdot\ldots\cdot a_r,\\
N_k<\frac{N_{k+1}}{a_{k+1}} \iff& \overline{b}_k\cdot a_k\cdot a_{k+1}\cdot\ldots\cdot a_r<
\overline{b}_{k+1}\cdot a_{k+2}\cdot\ldots\cdot a_r,\\
\iff&  \overline{b}_k\cdot a_k\cdot a_{k+1}<\overline{b}_{k+1},
\end{align*}
which is true since $\overline{b}_{k+1}=\overline{b}_k\cdot a_k\cdot a_{k+1}+b_{k+1}$ with $b_{k+1}>0$.
\end{proof}

We have the following result.

\begin{proposition}[{\cite[Prop 1.17]{PP:22}}]
Let $G$ be the semi-group generated by a Zariski sequence $(v_0,v_1,\dots v_r)$. Then the conductor $c$ of $G$ equals
\[
c=\sum_{k=1}^r (n_k-1)v_k-v_0+1.
\]
\end{proposition}

Let $f$ be an irreducible power series. Recall that the semi-group $\Gamma (f) $ of $f$ is defined~by
\[
\Gamma (f):=\{v_f(g), g \ \text{a power series such that } g \not\equiv 0 (f)\}.
\]

\begin{proposition} Assume $f$ is an irreducible power series, and $\mathcal{T}$ is its tree. The semi-group generated by $(v_0,\dots ,v_r)$ is the semi-group $\Gamma (f)$.
\end{proposition}

\begin{proof}
Let $g$ be an irreducible power series separating from $f$ on a vertex of the tree $\mathcal {T}$.  Then 
\[
i(f,g)=N_id_g=v_ia_id_g.
\]
Hence $i(f,g)\in \langle v_0,\dots ,v_r\rangle $.

Now assume that $g$ separates from $f$ on a dead end of the tree $\mathcal{T}$.  Then
\[
i(f,g)=d_g\cdot \overline{b}_i\cdot a_{i+1}\cdot\ldots\cdot a_r=d_g\cdot v_i.
\]
Hence $i(f,g)\in  \langle v_0,\dots ,v_r\rangle $.
Moreover, we see that if $d_g=1$, then $v_i\in \Gamma (f)$.

Finally, assume that $g$ separates from $f$ between the vertices $v_{i-1}$ and~$v_i$.
We have $i(f,g)>N_{i-1}$.

We want to show that $N_{i-1}>c_{i-1}\cdot a_i\cdot\ldots\cdot a_r$, where $c_{i-1}$ is the conductor of the semigroup generated by $\left\langle \frac{v_0}{a_i\cdot\ldots\cdot a_r}, \dots \frac{v_{i-1}}{a_i\cdot\ldots\cdot a_r}\right\rangle $:
\[
a_i\cdot\ldots\cdot a_r\cdot c_{i-1}=
N_{i-1}-\frac{N_{i-1}}{a_{i-1}}+N_{i-2}-\frac{N_{i-2}}{a_{i-2}}+N_{i-3}-\dots -\frac{N_2}{a_2}+N_1-v_0+1.
\]
Then $a_i\cdot\ldots\cdot a_r\cdot c_{i-1}<N_{i-1}$ and $i(f,g) >N_{i-1}>c_{i-1}\cdot a_i\cdot\ldots\cdot a_r$. 
Hence $i(f,g)\in \langle v_0,\dots v_{i-1}\rangle$.
\end{proof}

 \section{Multiplicity of the tree and \texorpdfstring{$\delta$}{delta}-invariant}\label{sec:delta}

In this section we want to compute the $\delta$-invariant of a series $f$ in terms of the multiplicity of its tree.
Let $f\in \bk[[x,y]]$ be a reduced power series.
Let $r(f)$ be the number of arrows not decorated by $(0)$ of the tree.
In this section, we shall use results from \cite{cnd} and \cite{cnv:14}. For this purpose, we need the following:

\begin{proposition}[{\cite[Prop 3.3]{cnv:14}}]
For all $v\in\cV$, we have
\[
N_v=\sum_{\alpha \in \arista\setminus \arista_0}\overline {\rho}(v,\alpha).
\]
\end{proposition}

This proves that the definitions of the multiplicity in \cite{cnd} and in the present work are the same, and we can use results of \cite{cnd}.

\begin{proposition}[\cite{cnd}]
The number $-M(\mathcal{T}_f)+r(f)$ is even.
\end{proposition}

Let
\[
\mathcal{O}_f:=\bk[[x,y]]/(f).
\]
Let $\overline{\mathcal{O}}_f$ be the integral closure of $\mathcal{O}_f$.
Define
\[
\delta (f):= \dim_\bk \overline{\mathcal{O}}_f/\mathcal{O}_f.
\]

We can state and prove the main result of the paper.
\begin{theorem}\label{thm:main}
$2\delta(f)=-M(f)+r(f)$.
\end{theorem}

\begin{proof}
Assume first that $f$ is an irreducible power series. Then $r(f)=1$.
We have
\[
M(\mathcal{T})=-\sum_{i=1}^r \left(N_i-\frac{N_i}{a_i}\right)-\frac{N_1}{a_1}.
\]
Then
$c(f)=-M(\mathcal{T})+1$.
It is proven in \cite[Thm~2.1]{PP:22} that
$c(f)=2\delta(f)$.
Then the result is true when $f$ is irreducible.

Let us recall \cite[Thm~2.1]{PP:22}:
\begin{quote}\em
If $f=f_1\dots f_r$ with irreducible coprime factors $f_i$, then
\[
\delta (f)=\sum_{i=1}^r \delta(f_i) +\sum _{1\leq i <j\leq r } i(f_i,f_j).
\]
\end{quote}

If we define
\[
\tilde{\delta}(\mathcal{T}):=\frac{-M(\mathcal{T})+r(f)}{2},
\]
then
\cite[Prop.~4.13]{cnd} says that
\[
\tilde{\delta}(\mathcal{T})=\sum _{\alpha \in \arista\setminus {\arista_0}}\tilde{\delta}(\mathcal{T}_{\alpha} )+
\frac{i(\arista\setminus {\arista_0},\arista\setminus {\arista_0})}{2},
\]
which proves the general case.
\end{proof}

\begin{remark}\label{rem:kou}
The combination of Theorems~\ref{thm:cv14} and~\ref{thm:main} 
provides a new interpretation
of Kouchnirenko's theorem, namely it refers to $\bar{\mu}$ and not $\mu$ for some
positive characteristics. Moreover, it can be extended without any hypothesis of non-degeneracy.
\end{remark}

\section{The multiplicity of the tree and the Milnor number}\label{sec:milnor}

Recall that the Milnor number of $f$ is
\[
\mu(f):=\dim_\bk \bk[[x,y]] / \left(\frac{\partial f}{\partial x}, \frac{\partial f}{\partial y}\right).
\]
If $\chr\bk=0$, Milnor results~\cite[Theorem~10.5]{milnor:68} imply that
\[
\mu(f) = 1- M(\cT)=2\delta(f)-r(f)+1,
\]
and Deligne~\cite{del:73} showed that in general we only have
\[
\mu(f)\geq 1 - M(f)=2\delta(f)-r(f)+1.
\]
Note that in positive characteristic, the Milnor numbers of $f$ and $uf$ ($u$~a unit) may differ, see e.g.
Examples~\ref{ejm:milnor1} and~\ref{ejm:milnor2}, but $\overline{\mu}$ does not  \cite[Prop.~ 3.1]{PP:22}.
For a review about these equalities (and inequalities) we refer to~\cite{evpl:18}.

\color{black}

\begin{cjt}
Let $\cT$ be a minimal tree of $f$.
The equality holds if and only if $\chr\bk$ divides $N_v$ for no $v\in\cV\cup\cA_0$.
\end{cjt}

We shall now show that the conjecture was already proven in some particular cases. 

\begin{definition}
Let $\Delta$ be a face of the Newton polygon of $f$. We call $f$ \emph{non-degenerate} (ND) along $\Delta$ if the Jacobian ideal of $f_{\Delta}$ has no zero in the torus $(\bk^*)^2$.
We say that $f$ is \emph{Newton non-degenerate} (NND) if $f$ is ND along each face (of any dimension ) of the Newton polygon of $f$. 
\end{definition}

In 1976, Kouchnirenko~\cite{kou:76}  proved that if $f$ is NND and convenient, then $\mu(f)=1-M(f)$.  This result was extended by~\cite{bgm:12} in 2010 (published in~2012) without the hypothesis of convenient (if $\mu(f)<\infty$.).

\begin{proposition}\label{prop:nnd}
Let $f$ be non-degenerate. If $\chr \bk$ divides $N_v$ for no $v$, then $f$ is NND.
\end{proposition}

\begin{proof}
Let $\Delta$ be a face of dimension one of the Newton polygon of $f$.  Let $v$ be the corresponding vertex on the tree, and $N_v$ the corresponding multiplicity. The equation of the face  is $aX+ bY=N_v$. 
We want to show  that if $\chr \bk=p$ does not divide  $N_v$, then $f$ is \nd{} along $\Delta$. We have
\begin{align*}
f_{\Delta}(x,y)=&x^ny^m\sum_i ^k c_i x^{ia}y^{(k-i)b},\\
\dfrac{\partial f_{\Delta}}{\partial x}(x,y)=&\sum_i^k c_i (n+ia)x^{n+ia-1}y^{(k-i)b+m},\\
\dfrac{\partial f_{\Delta}}{\partial y}(x,y)=&\sum _i^k c_i (m+(k-i)b)x^{n+ia}y^{(k-i)b+m-1}.
\end{align*}
If the Jacobian ideal of $f_{\Delta}$ has a zero $(\alpha, \beta)$ in $(\bk^*)^2$, we have 
\begin{align*}
\sum_i^k c_i (n+ia)\alpha^{ia}\beta^{(k-i)b}&=0,\\
\sum _i^k c_i (m+(k-i)b)\alpha^{ia}\beta^{(k-i)b}&=0,
\end{align*}
which is impossible: for all $i$ we have 
\[
(n+ia)(m+(k-i+1)b-(n+(i-1)a)(m+(k-i)b))=N_v\neq 0,
\]
since $p$ does not divide $N_v$.

Now if $f$ is not ND with respect to a vertex of the Newton polygon, then $p$ divides $N_v$ for a face of the Newton polygon which contains the vertex. 
\end{proof}

\begin{corollary}
If $f$ is non-degenerate and $\chr\bk$ divides $N_v$ for no $v$, then $\mu(f) =2\delta(f)-r(f)+1$.
\end{corollary}

\begin{proof}
Assume $f$ is non-degenerate.  Recall that we assume that the tree is minimal.
If $f$ is non-degenerate and 
$\chr\bk$ divides $N_v$ for no $v$,
 then $f$ is NND (Proposition~\ref{prop:nnd}) and $\mu (f)= 1-M(f)$~\cite{bgm:12}.
\end{proof}

\begin{corollary}
The conjecture is true if $f$ is non-degenerate.
\end{corollary}

\begin{proof}
If $p$ divides some $N_v$, then $f$ is not \nd{} with respect of $\Delta$, the corresponding face of the Newton polygon. Using {\cite[Prop~2.12]{gn:12}} we have that $\mu(f)>2\delta(f)-r(f)+1$.
\end{proof}

Now we recall two interesting results from \cite{evpl:18}. Let $p:=\chr\bk$, $l(x,y):=ax+by$, and $P_l(f):=b\frac{\partial f}{\partial x}-a\frac{\partial f}{\partial y}$.

\begin{proposition}[{\cite[Prop. 1.1]{evpl:18}}]
Let $l$ be a regular parameter of $\bk[[x,y]]$, and $f=f_1\cdot f_2\cdot\ldots\cdot f_r$ such that $i(f_i,l)\not\equiv 0 \bmod p$. Then 
\[
i(f,P_l(f))=2\delta(f)+i(f,l)-r=-M(\mathcal{T}(f))+i(f,l).
\]
\end{proposition}

\begin{proposition}[{\cite[Prop. 3.4]{evpl:18}}]
Assume that there exists a regular parameter such that  $i(f,l)=\ord (f)$ and $i(f,P_l(f))<p$. Then 
$\mu(f)=-M(\mathcal{T}(f))+1=2\delta-r+1$. 
\end{proposition}

Then we deduce the following:
\begin{proposition}\label{prop:77}
If $p>-M(\mathcal{T}(f))+\ord(f)$, then $\mu(f)=-M(\mathcal{T})+1=2\delta(f) -r+1$.
\end{proposition}

Now we study the case where $f$ is irreducible.
In this case, we have two results. Garc{\'i}a-Barroso and P{\l}oski~\cite[Theorem~4.1]{evpl:18} proved  the following:
\begin{quote}\em
Let $n^*:=\max (a_1,\dots,a_r)$ (See Figure{\rm~\ref{fig:branche})}. If  $\chr\bk=p>n^*$, then 
\[
\mu(f)=2\delta(f)\Longleftrightarrow\forall k, 0\leq k\leq r, \nu_k \not\equiv 0 \bmod p.
\]
\end{quote}
This proves the conjecture when $p>n^*$.
On the other hand Hefez, Rodrigues, Solom\~{a}o proved in~\cite{hrs:18} that if 
$\forall k\in\{0,\dots,r\}$ we have that $\nu_k \not\equiv 0 \bmod p, $ then $\mu(f)=2\delta(f)$.

\begin{example}
We consider 
\[
f(x,y):=(x-a_1y)(x-a_2y)(x-a_3y)(x-a_4y)+xy^5+x^4y.
\]
Let
\[
b_i:=a_{i+1}-a_1,\ 1\leq i<4,\quad 
c_i:=a_{i+2}-a_2,\ 1\leq i<3,\quad 
d_1:=a_4-a_3.
\]
We assume $p\neq 2$.

\begin{enumerate}[label=\rm(\arabic{enumi}), series=tt, leftmargin=*, widest=14]

\item We assume all the $a_i$ are pairwise distinct $\!\!{}\bmod{p}$.
In this case, the tree is $\mathcal{T}_1$ in Figure~\ref{subfig:t1}, 
and we have 
\[
\mu(f)=1-M(\mathcal{T}_1)=2\delta(f)-3=9,\quad\forall p\neq 2.
\]

\item We assume that $a_1\equiv a_2\equiv 0\bmod{p}$, $a_3,a_4\not\equiv 0\bmod{p}$, $a_3\not\equiv a_4\bmod{p}$.
In this case, the tree is  $\mathcal{T}_2$ in Figure~\ref{subfig:t2},
and we have 
\[
\mu(f)=1-M(\mathcal{T}_2)=2\delta(f)-3=13,\quad\forall p\neq 2.
\]

\item Only $a_4$ is non-vanishing $\!\!{}\bmod{p}$.
In this case, the tree is $\mathcal{T}_3$ in Figure~\ref{subfig:t3},
and we have 
\[
\mu(f)=1-M(\mathcal{T}_3)=2\delta(f)-3=15,\quad\forall p\neq 2,7,\quad 
\mu=17\text{ if } p=7.
\]

\item We assume that $a_i\equiv 0\bmod{p}$.
In this case, the tree is $\mathcal{T}_4$ in Figure~\ref{subfig:t4},
and we have 
\[
\mu(f)=1-M(\mathcal{T}_4)=2\delta(f)-1=17,\quad\forall p\neq 2,5,\quad 
\mu=20\text{ if } p=5.
\]

\begin{figure}[ht]
\centering
\begin{subfigure}[b]{.24\textwidth}
\centering
\begin{tikzpicture}
\foreach \x in {-30, -10, 10, 30}
{
\draw[->]  (0,0) -- (\x:1);
}
\filldraw[fill=white] (0,0) node[left] {(4)} circle [radius=.1cm];
\end{tikzpicture}
 \subcaption{$\mathcal{T}_1$}
\label{subfig:t1}
\end{subfigure}
\begin{subfigure}[b]{.24\textwidth}
\centering
\begin{tikzpicture}
\draw[->] (0,0) -- (0, 2);
\draw[->] (0,0) -- (-20:1);
\draw[->] (0,0) -- (20:1);
\draw[->] (0,1.25) -- (1,1.25);
\filldraw[fill=white] (0,0) node[left=5pt] {(4)} node[above left] {$1$} circle [radius=.1cm];
\filldraw[fill=white] (0,1.25) node[left=5pt] {(8)} node[above left] {$1$}node[below left] {$3$} circle [radius=.1cm];
\end{tikzpicture}
 \subcaption{$\mathcal{T}_2$}
\label{subfig:t2}
\end{subfigure}
\begin{subfigure}[b]{.24\textwidth}
\centering
\begin{tikzpicture}
\draw[<->] (0,-1) -- (0, 1);
\draw[->] (0,0) -- (-20:1);
\draw[->] (0,0) -- (20:1);

\filldraw[fill=white] (0,0) node[left=5pt] {(7)} node[above left] {$1$}node[below left] {$2$} circle [radius=.1cm];

\end{tikzpicture}
 \subcaption{$\mathcal{T}_3$}
\label{subfig:t3}
\end{subfigure}
\begin{subfigure}[b]{.24\textwidth}
\centering
\begin{tikzpicture}
\draw[<->] (-1,0) node[left] {$(0)$} -- (1,0);
\draw[->] (0,0) -- (0,-1);

\filldraw[fill=white] (0,0) node[above right] {(20)} node[above left] {$5$}node[below right] {$3$} circle [radius=.1cm];

\end{tikzpicture}
 \subcaption{$\mathcal{T}_4$}
\label{subfig:t4}
\end{subfigure}
\caption{\footnotesize}

\end{figure}

\end{enumerate}
In all these cases, we use the fact that $f$ is non-degenerate.
\begin{enumerate}[label=\rm(\arabic{enumi}), resume=tt, leftmargin=*, widest=14]

\item No $a_i$ vanishes $\!\!{}\bmod{p}$, $a_1\equiv a_2\bmod{p}$, and $a_2,a_3,a_4$
are pairwise distinct $\!\!{}\bmod{p}$.
In this case the tree is $\mathcal{T}_{1,2}$ in Figure~\ref{subfig:t12},
and we have 
\[
\mu(f)=1-M(\mathcal{T}_{1,2})=2\delta(f)-2=10,\text{ if } p=0,\quad 
\mu=11\text{ if } p=5.
\]
From Proposition~\ref{prop:77}, we know that if $p> 13$, then  $\mu =10$.
This value is also obtained for the remaining cases $p=3,7,11,13$.

\item No $a_i$ vanishes $\!\!{}\bmod{p}$, $a_1\equiv a_2\not\equiv a_3\equiv a_4\bmod{p}$.
In this case the tree is $\mathcal{T}_{1,3}$ in Figure~\ref{subfig:t13},
and we have 
\[
\mu(f)=1-M(\mathcal{T}_{1,3})=2\delta(f)-1=11,\text{ if } p=0,\quad 
\mu=13\text{ if } p=5.
\]
Again, using  Proposition~\ref{prop:77}, we know that if $p> 14$, then  $\mu =11$, and we verify for $p=3,7,11,13$ that $\mu =11$ (using \texttt{Singular} or \texttt{Sagemath}).

\begin{figure}[ht]
\centering
\begin{subfigure}[b]{.32\textwidth}
\centering
\begin{tikzpicture}
\draw[->] (0,0) -- (-30:1);
\draw[->] (0,0) -- (0:1);
\draw[->] (0,0) -- (30:1.75);
\draw[->] (30:1) -- ($(30:1)+(-60:1)$) node[right] {$(0)$};
\filldraw[fill=white] (0,0) node[left] {$(4)$} circle [radius=.1cm];
\filldraw[fill=white] (30:1) node[above=2pt] {$(10)$}
node[left=3pt] {$3$} circle [radius=.1cm];
\node at ($(30:1)+(.3,-.1)$) {$2$};
\end{tikzpicture}
 \subcaption{\footnotesize$\mathcal{T}_{1,2}$}
\label{subfig:t12}
\end{subfigure}
\begin{subfigure}[b]{.32\textwidth}
\centering
\begin{tikzpicture}
\draw[->] (0,0) -- (-45:1.75);

\draw[->] (0,0) -- (45:1.75);
\draw[->] (45:1) -- ($(45:1)+(-45:1)$) node[right] {$(0)$};
\draw[->] (-45:1) -- ($(-45:1)+(-135:1)$) node[left] {$(0)$};
\filldraw[fill=white] (0,0) node[left] {$(4)$} circle [radius=.1cm];
\filldraw[fill=white] (45:1) node[above=2pt] {$(10)$}
node[left=3pt] {3} circle [radius=.1cm];
\filldraw[fill=white] (-45:1) node[right=2pt] {$(10)$} circle [radius=.1cm];
\node at ($(45:1)+(.4,-.1)$) {$2$};
\node at ($(-45:1)+(-.1,-.35)$) {$2$};
\node at ($(-45:1)+(-.1,.35)$) {$3$};
\end{tikzpicture}
 \subcaption{\footnotesize$\mathcal{T}_{1,3}$}
\label{subfig:t13}
\end{subfigure}
\begin{subfigure}[b]{.32\textwidth}
\centering
\begin{tikzpicture}
\draw[<->] (-1,0)  -- (1,0);
\draw[->] (0,0) -- (0,-1) node[below] {$(0)$} ;

\filldraw[fill=white] (0,0) node[above right] {$(15)$} node[above left] {$4$}node[below left] {$3$} circle [radius=.1cm];

\end{tikzpicture}
 \subcaption{\footnotesize$\mathcal{T}_{1,4}$}
\label{subfig:t14}
\end{subfigure}

\caption{\footnotesize}
\end{figure}

\item No $a_i$ vanishes $\!\!{}\bmod{p}$, $a_1\equiv a_2\equiv a_3\not\equiv a_4\bmod{p}$.
In this case the tree is $\mathcal{T}_{1,4}$ in Figure~\ref{subfig:t14},
and we have 
$\mu(f)=1-M(\mathcal{T}_{1,4})=2\delta(f)-1=11$ if $p=0$,
$\mu=12$ if $p=5$, and
$\mu=13$ if $p=3$.
We know that if $p>14$, $\mu =11$, and we verify that it is also the case for $p=7,11,13$.

\item No $a_i$ vanishes $\!\!{}\bmod{p}$, and they are equal $\!\!{}\bmod{p}$.
In this case the tree is $\mathcal{T}_{1,5}$ in Figure~\ref{subfig:t15},
and we have 
\[
\mu(f)=1-M(\mathcal{T}_{1,5})=2\delta(f)-1=12,\text{ if } p=0,\quad 
\mu=13\text{ if } p=5.
\]
We know that if $p>15$, then $\mu =12$, and  we verify that it is also the case for $p=3,7,11,13$.

\item We assume that  $a_1\equiv a_2\equiv 0\bmod{p}$ and $a_3\equiv a_4\not\equiv 0\bmod{p}$.
In this case the tree is $\mathcal{T}_{2,1}$ in Figure~\ref{subfig:t21},
and we have 
\[
\mu(f)=1-M(\mathcal{T}_{2,1})=2\delta(f)-2=14,\text{ if } p=0,\quad 
\mu=15\text{ if } p=5.
\]
We know  that  for $p>17$, $\mu =14$, and  we verify that it is also the case for $p=3,7,11,13,17$.

The conjecture is true for this family.

\begin{figure}[ht]
\centering
\begin{subfigure}[b]{.32\textwidth}
\centering
\begin{tikzpicture}
\draw[<->] (-1,0) node[left] {$(0)$} -- (1,0);
\draw[->] (0,0) -- (0,-1) node[below] {$(0)$} ;

\filldraw[fill=white] (0,0)  node[above left] {$5$}node[below left] {$4$} circle [radius=.1cm];

\end{tikzpicture}
 \subcaption{\footnotesize$\mathcal{T}_{1,5}$}
\label{subfig:t15}
\end{subfigure}
\begin{subfigure}[b]{.32\textwidth}
\centering
﻿\begin{tikzpicture}
\draw[->] (0,0) -- (0, 2);
\draw[->] (0,0) -- (0,-1) node[below] {$(0)$};
\draw[->] (0,0) -- (1,0);
\draw[->] (0,1.25) -- (1,1.25);
\filldraw[fill=white] (0,0) node[left=5pt] {$(10)$} node[above left] {$3$} node[below right] {$2$}  circle [radius=.1cm];
\filldraw[fill=white] (0,1.25) node[left=5pt] {$(8)$} node[above left] {$1$} node[below left] {$3$} circle [radius=.1cm];
\end{tikzpicture}
 \subcaption{\footnotesize$\mathcal{T}_{2,1}$}
\label{subfig:t21}
\end{subfigure}
\begin{subfigure}[b]{.32\textwidth}
\centering
\begin{tikzpicture}
\draw[->] (0,0)  node[above left] {$1$} node[below left] {$2$} node[left=5pt] {$(7)$} -- (2,0);
\draw[<->] (0,-1) -- (0, 1) ;
\draw[->] (1,0) node[above left] {$5$} node[below right] {$2$} node[above right] {$(16)$}  -- (1, -1) node[below] {$(0)$} ;
\filldraw[fill=white] (0,0) circle [radius=.1];
\filldraw[fill=white] (1,0) circle [radius=.1];
\end{tikzpicture}
 \subcaption{\footnotesize$\mathcal{T}_{1,6}$}
\label{subfig:t16}
\end{subfigure}
\caption{\footnotesize}
\end{figure}

\end{enumerate}

Now we consider the case where $p=2$.

\begin{enumerate}[label=\rm(\arabic{enumi}), resume=tt, leftmargin=*, widest=14]
\item We assume that $a_i\neq 0 $ for all $i$. Then the tree is $\mathcal{T}_{1,5}$. The multiplicity of the tree is $11$, and we can compute $\mu=20$.
\item We assume that exactly one of the $a_i=0$. Then the tree is $\mathcal{T}_{1,4}$ we have $M(\mathcal{T}_{1,4})=10$ and $\mu=11$ since $2$ does not divide $N_v$.
\item We assume that exactly two of the $a_i$ vanish. The tree is $\mathcal{T}_{2,1}$. Its multiplicity is $13$ and $\mu=20$.
\item We assume that exactly three of the $a_i$ vanish. The tree is $\mathcal{T}_{1,6}$ in Figure~\ref{subfig:t16}.
Its multiplicity is $15$ and we compute $\mu=19$.
\item We assume all the $a_i=0$. The tree is $\mathcal{T}_{4}$. The multiplicity is $16$ and $\mu=20$.
\end{enumerate}

\end{example}

\begin{example}\label{ejm:milnor1}
Let us consider Example~\ref{ex1}.
In this example, $f$ is irreducible. The results of \cite{evpl:18} tell us that if $p>3$ then $\mu=156$ if and only if $p\neq 5,101$, where $\mu=157$.  We can check that for $p=2$, $\mu=\infty$, and $p=3$, $\mu=166$. Note also that if we multiply $f$ by a \emph{random} unit then $\mu=168$ for $p=2$ and $\mu=157$ for $p=3$; nothing changes for the other primes. Then the conjecture is true for this example.
\end{example}

\begin{example}\label{ejm:milnor2}
Let us consider Example~\ref{ex2}.
First we assume $p\neq 2$. If $p>111$, then $\mu=102$. The prime numbers which divide $N_v$ for some $v$ are $p=7$, and we have $\mu=105$, $p=11$ with $\mu=\infty$ and $p=13$ with $\mu=104$. 
For the remaining primes $p<111$, we also have  $\mu=102$.
For $p=2$, we have $\mu=133$. Note also that if we multiply $f$ by a some unit then $\mu=118$ for $p=2$, $\mu=104$ for $p=3$, and $\mu=105$ for $p=11$; nothing changes for the other primes. 
Then the conjecture is true in this example.
\end{example}
 
\subsection*{Acknowledgments}
First named author is partially supported by MCIN/AEI/ 10.13039/501100011033 (grant code: PID2020-114750GB-C31),
and by Departamento de Ciencia, Universidad y Sociedad del Conocimiento del Gobierno de Arag{\'o}n
(grant code: E22\_20R: ``{\'A}lgebra y Geometr{\'i}a'').
Second named author is partially supported by PID2020-114750GB-C31
and PID2020-114750GB-C32.

% \bibliographystyle{amsalpha}
% \bibliography{biblio}

\providecommand{\bysame}{\leavevmode\hbox to3em{\hrulefill}\thinspace}
\providecommand{\MR}{\relax\ifhmode\unskip\space\fi MR }
\providecommand{\MRhref}[2]{%
  \href{http://www.ams.org/mathscinet-getitem?mr=#1}{#2}
}
\providecommand{\href}[2]{#2}

\end{document}